\documentclass[a4paper, 10pt]{article}
\usepackage{latexsym,amsthm,amssymb,amscd}
\usepackage{color,graphicx}
\usepackage[numbers,square,sort&compress]{natbib}
\author{{{\bf M. Mohammadnezhad, S. Golalizadeh, M. Boostan,}} {\bf N. Soltankhah \thanks{Corresponding author: soltan@alzahra.ac.ir, soltankhah.n@gmail.com }} \\
  {\footnotesize {\bf Department of Mathematics, Faculty of Mathematical Sciences, Alzahra University, Tehran, Iran } }}
\title {Super-simple $(v,5,2)$ directed designs and  their smallest defining sets with its application in LDPC codes  }
\date{}
\newtheorem{theorem}{Theorem}
\newtheorem{definition}{Definition}
\newtheorem{proposition}{Proposition}

\newtheorem{lemma}{Lemma}

\newtheorem{example}{Example}
\newtheorem{construction}{Construction}

\newcommand{\B}{\mathcal{B}}
\newcommand{\G}{\mathcal{G}}
\newcommand{\D}{\mathcal{D}}
\usepackage{ptext}
\usepackage{subcaption}
\usepackage{tikz,amsmath} 
\tikzstyle{vertex}=[circle, draw, inner sep=.5pt, minimum size=0pt] 
\newcommand{\vertex}{\node[vertex]}
\begin{document}
%\baselineskip=28
\maketitle 
\begin{abstract}
In this paper, we show that for all $v\equiv 0,1$ (mod 5) and $v\geq 15$, there exists
a super-simple $(v,5,2)$ directed design, also for these parameters there exists a super-simple $(v,5,2)$ directed design such
that its smallest defining sets contain at least half of its blocks. Also, we show that these designs are useful in constructing parity-check matrices of LDPC codes.
\end{abstract}
\section{Introduction and preliminaries}
A group divisible design ( or GDD) is a triple $(X,\mathcal{G},\mathcal{B})$ which satisfies the following properties:
\begin{enumerate}
\item
$\mathcal{G}$ is a partition of a set $X$ into subsets called groups;
\item
$\mathcal{B}$ is a set of subsets of $X$ called blocks such that a group and a block intersect in at most one point;
\item
each pair of points from distinct groups occurs in exactly $\lambda$ blocks.
\end{enumerate}

The group type of GDD is the multiset $\{|G|: G\in \mathcal{G}\}$. We use the notation $g_1^{u_1}g_2^{u_2}\cdots g_n^{u_n}$ to denote $u_i$ occurrences of $g_i$ for $1\leq i\leq n$ in the multiset. A GDD with block sizes from a set
of positive integers $K$ is called a $(K,\lambda)$-GDD. When $K=\{k\}$, we simply write $(k,\lambda)$-GDD. When $\lambda=1$, we simply write $K$-GDD. A $(K,\lambda)$-GDD with group type $1^v$ is called a pairwise balanced design and denoted by PBD$(v,K,\lambda)$. A $(k,\lambda)$-GDD with group type $1^v$ is called a balanced incomplete block design, denoted by $(v,k,\lambda)$-BIBD.

Some generalizations have been introduced for the concept of designs. Gronau and
Mullin \cite{gronau} for the first time, introduced a new definition of block
designs called super-simple block designs. A super-simple $(v,k,\lambda)$  design  is a block design such that any two blocks of the design
intersect in at most two points. A simple block design is a block design such that it has no repeated blocks.
The existence of super-simple $(v,4,\lambda)$ designs have been
characterized for $2\leq \lambda \leq 9$ 
%except $\lambda =7$
, see \cite{chen3, gchen, chen6, chen1, 
chen5, chen0, gronau,sun, zhang}. Also, the existence of super-simple
$(v,5,\lambda)$ designs have been characterized for $2\leq \lambda \leq 5$,   see  \cite{chen11, chen2, chen4, hans}.

A directed group divisible design $(K,\lambda)$-DGDD is a group divisible design
in which every block is ordered and each ordered pair formed from distinct
elements of different groups occurs in exactly $\lambda$ blocks. A $(k,\lambda)$-DGDD with group type $1^v$ is called a directed balanced incomplete block design and denoted by $(v,k,\lambda)$-DBIBD or $(v,k,\lambda)$DD.
 A $(K,\lambda)$-DGDD is super-simple if its underlying $(K,2\lambda)$-GDD is super-simple.\\
A transversal design, TD$(k,\lambda,n)$, is a $(k,\lambda)$-GDD of group type $n^k$. When $\lambda=1$, we use the notation TD$(k,n)$. 
\begin{lemma}\label{l11111}\cite{abel}
\begin{enumerate}
\item[1]
A TD$(q+1,q)$ exists, consequently, a TD$(k, q)$ exists for any positive integer $k (k\leq q + 1)$, where $q$ is a prime power.
\item[2]
A TD$(7,n)$ exists for all $n\geq 63$.
\end{enumerate}
\end{lemma}
A set of blocks which is a subset of a unique $(v,k,\lambda)$DD is said to be a 
defining set of the directed design. A minimal defining set is a defining set, no
proper subset of which is a defining set. A smallest defining set, is a defining
set with the smallest cardinality.
A $(v,k,t)$ directed trade of volume $s$ consists of two disjoint collections $T_1$ and $T_2$ each of $s$ ordered $k$-tuples of a $v$-set $X$ called blocks, such that every ordered $t$-tuple of distinct elements of $X$ is covered by exactly the same number of blocks of $T_1$ as of $T_2$. Such a directed trade is
usually denoted by $T=T_1-T_2$. In a $(v,k,t)$ directed trade, both collections of blocks cover the same set of elements. This set of elements is called the foundation of the trade. In \cite{soltan}, it has been shown that the minimum volume of a $(v, k, t)$ directed trade is $2^{\lfloor{\frac{t}{2}\rfloor}}$ and that directed trades
of minimum volume and minimum foundation exist.
In this paper we use a special type of a directed trade which is defined as follows. 
\begin{definition}\cite{soltankhah1}
Let $T=T_1-T_2$ be a $(v,k,2)$ directed trade of volume $s$ with blocks $b_0$, $b_1$,$\cdots$, $b_{s-1}$ such that each pair of consecutive blocks of $T_1$ ($b_i$, $b_{i+1}$ $i=0,1,\cdots,s-1$ (mod $s$)) is a trade of volume $2$. Such a trade is called a cyclical trade. (We denote a cyclical trade of volume $ s $ by $ CT_s $)
\end{definition}
%\begin{definition}
%A cyclical trade of volume $s$ is a set of $s$ blocks $ \lbrace b_0, b_1, \cdots, b_{s-1}\rbrace $ if each pair of consecutive blocks $ b_i,b_{i+1} $ for $ i=0,1,\cdots,s-1 $, have the same pair that forms a $ (v,k,t) $ directed trade $ T=T_1-T_2 $ of volume 2, and in total $ s $ directed trades of volume 2. (We denote a cyclical trade of volume $ s $ by $ CT_s $) 
%\end{definition}
If $\mathcal{D}=(V,\mathcal{B})$ is a directed design and if $T_1\subset \mathcal{B}$, we say that $\mathcal{D}$ contains the directed trade $T$. Defining sets for directed designs are strongly related to trades. This relation is illustrated by the following result.
\begin{proposition}
Let $\mathcal{D}=(V,\mathcal{B})$ be a $(v,k,\lambda)$DD and let $S\subset \mathcal{B}$, then $S$ is a defining set of $\mathcal{D}$ if and only if $S$ contains a block of every $(v,k,2)$ directed trade $T=T_1-T_2$ such that $T$ is contained in $\mathcal{D}$.
\end{proposition}
Each defining set of a $(v,k,\lambda)$DD $\mathcal{D}$, contains at least one block from each trade in $\mathcal{D}$. In particular, if $\mathcal{D}$ contains $m$
mutually disjoint directed trades then the smallest defining set of $\mathcal{D}$ must contain at least $m$ blocks. If a directed design $\mathcal{D}$ contains a cyclical trade  of volume $s$, then each defining set for $\mathcal{D}$ must contain at least $\lfloor\frac{s+1}{2}\rfloor$ blocks of $T_1$.

Some results have been obtained on  $(v,k,\lambda)$DDs for special $k$ and $\lambda$  and their defining sets. For example,  in \cite{es.}, it has been proved that if $\D$ is a $(v,3,1)$DD, then a defining set of $\D$ has at least $\frac{v}{2}$ blocks. In \cite{Grannell}, it has been shown that for each admissible value of $v$, there exists  a simple $(v,3,1)$DD whose smallest defining sets have at least a half of the blocks. In
\cite{soltankhah1}, it has been shown that the necessary and
sufficient condition for the existence of a super-simple
$(v,4,1)$DD is $v\equiv 1$ (mod 3) and for these values of $v$ except $v=7$, there exists a super-simple $(v,4,1)$DD whose  smallest
defining sets have at least a half of the blocks. Also,
in \cite{soltankhah2}, it has been shown  that for all
$v\equiv 1,5$ (mod 10) except $v=5,15$, there exists a super-simple
$(v,5,1)$DD such that their smallest defining sets have at least a
half of the blocks. In \cite{goli}, the authors showed   that for all
$v\equiv 1$ (mod 3), there exists a super-simple
$(v,4,2)$DD such that their smallest defining sets have at least a
half of the blocks.
\\
In this paper, we prove that the necessary and sufficient condition
for the existence of a super-simple $(v,5,2)$DD is $v\equiv 0,1$ (mod
$5$) $(v\geq 15)$ and for these values of $v$, there exists  a super-simple
$(v,5,2)$DD  whose  smallest defining sets have at least a
half of the blocks. We introduce the following quantity
$$d=\frac{the\  total \  number \ of \  blocks \   in \   a \   smallest  \   defining \  set \   in \   \D}{the \  total \   number \   of \  blocks \  in \  \D}$$
and we show for all admissible values of $v$, $d\geq \frac{1}{2}$.\\
In the last section we provide a new method to construct parity-check matrices of LDPC codes by these  designs.
%2. Define $\lambda_i$($\rho_i$) to be the fraction of message (check) edges of degree i and let
%$$\lambda(x) = \Sigma_i \lambda_i x^{i-1} (\rho(x) = \Sigma_j \rho_j x_{j-1})$$.
%3. Techniques for the design of good irregular LDPC codes by finding good sequences $ \lbrace\lambda_i,\rho_i\rbrace $ are given in [1898] for the BEC and in [1487] for more general
%channels.
\section{Recursive Constructions}
For some values of $v$, the existence of a super-simple $(v,5,2)$DD will be proved by recursive constructions that which are presented  in this section for later use.
\begin{construction}(Weighting)\label{1}\cite{goli}
 Let $ (X,\G,\B) $ be a super-simple DGDD with index $\lambda_1$ and with $d\geq \frac{1}{2}$. Let $ w:X\rightarrow Z^+ \bigcup \{0\} $ be a weight function on $X$, where $Z^+$ is the set of positive integers. Suppose that for each block $ B\in \B $, there exists a super-simple $ (k,\lambda_2)$-DGDD  of type $ \{w(x): x\in B\} $ with $d\geq \frac{1}{2}$. Then there exists a super-simple $(k,\lambda_1\lambda_2)$-DGDD of type $\{\sum_{x\in G_i} w(x):  G_i\in \G \}$ with $d\geq \frac{1}{2}$.
\end{construction}

\begin{construction}\label{2}\cite{goli}
If there exist a super-simple $(k, \lambda)$-DGDD of type $g_1^{u_1}\cdots g_t^{u_t}$ with $d\geq \frac{1}{2}$
 and a super-simple
$(g_i +\eta, k,\lambda)$DD for each $ i (1\leq i\leq t)$ with $d\geq \frac{1}{2}$, then there exists a super-simple $(\sum_{i=1}^t  g_iu_i+\eta, k, \lambda)$DD with $d\geq \frac{1}{2}$, where
$\eta = 0 \ \  or  \ \ 1$.
\end{construction}
%\begin{construction}\label{3}\cite{goli}
%If there exists a $K-GDD $ of type $g_1^{u_1}g_2^{u_2}...g_n^{u_n}$, a super-simple $(\alpha g_i +1,5,2)$DD for each i, $i=1,2,...,n$ and a super-simple $4$-DGDD of type $\alpha^k$ for each $k \in K$, then there exists a super-simple $(\alpha \sum_{i=1} ^n g_i u_i +1,5,2)$DD.
%\end{construction}
\section{Direct Construction}
In this section, we construct some super-simple $(v,5,2)$DDs for some small admissible values of $v$ and some super-simple directed group divisible designs by direct construction and for these values of $v$, we show that the parameter $d$ for constructed designs is at least $\frac{1}{2}$. \\
In what follows we use the notation $+d$ (mod $v$), which denotes that all elements of the base blocks should be
developed cyclically by adding $d$ (mod $v$) to them, while the infinite point $\infty$, if it occurs in the base blocks, is always
fixed.We usually omit $+d$ when $d = 1$.\\
Let $[a,b]^{0,1}_5$ be the set of positive integers $v$ such that $v\equiv 0,1$ (mod $5$) and $a\leq v\leq b$.
\begin{lemma}\label{l1} There exists a super-simple $(v,5,2)$DD for all $v\in[15,86]^{0,1}_5\cup \{95,110,111,115,116,130,131\} $,
whose  smallest defining sets have at least a half of the
blocks.
\begin{proof}
For $ v=15$ and $\mathcal{G}=Z_{14}\cup\lbrace \infty \rbrace $, The following base blocks by +2 (mod 14)  form a super-simple $(15,5,2)$DD.
\begin{center}
\begin{tabular}{ccc}
  (1,0,2,3,8) & (0,3,13,11,9) & (0,1,4,10,9)  \\
  (0,7,11,4,2) & (1,0,$\infty$,5,7) & (13,2,$\infty$,0,10)  \\
\end{tabular}
\end{center}
This design contains 42 blocks, each of  three  columns  has 7 
disjoint directed trades of volume 2. Since each defining set for this design
 must contain one 5-tuple of  each directed trade in each of columns, then each
 defining set contains at least $7\times 3=21$ blocks. So $d\geq \frac{1}{2}$.\\
 For $v=25$ and $\mathcal{G}=Z_{24}\cup\lbrace \infty \rbrace $, the following base blocks by $+1$ (mod $24$) form a super-simple $(25,5,2)$DD.
 \begin{center}
\begin{tabular}{ccc}
  (0,5,1,7,15) & (22,0,5,21,11) & (12,0,1,10,4)  \\
  (2,0,$\infty$ ,17,21) & (13,6,1,0,9) &  \\
\end{tabular}
\end{center}
 There are $120$ blocks in a super-simple $(25,5,2)$DD. The first two columns have $48$ disjoint directed trades of volume $2$, and
the last column is a cyclical trade of volume $24$. Since each defining set for this super-simple directed design must contain at
least one $5$-tuple of each directed trade in the first two columns and $12$ $5$-tuples of cyclical trade in the last column, then
each defining set must contain at least $48+12=60$ blocks. Therefore for this super-simple $(25,5, 2)$DD the inequality $d\geq \frac{1}{2}$ is satisfied.\\
 For $v\in [16,36]^{0,1}_5$ except $v=15,25$, the results are summarized in the following table. 
\begin{center}
\begin{tabular}{|c|cccc|c|c|c|}
\hline
$v$ & base blocks & & & & mod & $b_v$ & $d\geq$ \\
  \hline
  16 & (3,0,1,8,6) & (1,7,14,2,11) & (2,5,0,1,4) &  & +2 mod $16$ & $48$ & $\frac{3\times 8}{48}$\\
    & (7,0,3,11,5) & (0,2,12,7,8) & (1,0,10,7,9)  &  &  & & \\
    \hline
 20 &  (0,4,3,9,16) & (9,0,1,18,14) & & & mod $19$ & 76 & $\frac{2\times 19}{76}$\\
 & (5,0,$\infty$, 7,8) & (11,2,4,8,0) & & & & &\\ 
 \hline
21 & (3,0,6,8,7) & (0,19,18,7,10) & & & mod $21$ & 84 & $\frac{2\times 21}{84}$\\
 & (8,2,4,0,16) & (0,9,14,4,15) & & & & &\\ 
 \hline
%25 & (0,5,1,7,15) & (22,0,5,21,11) & (12,0,1,10,4) & & mod $24$ & 120 & $\frac{2\times 24+24}{120}$\\
%&  (2,0,$\infty$ ,17,21) & (13,6,1,0,9) & & & & &\\ 
 % \hline
26 & (5,13,0,7,22) & (3,0,19,7,25) & (8,0,13,14,24) & & mod $26$ & 130 & $\frac{2\times 26+13}{130}$\\
 & (16,0,11,20,19) & (0,12,1,24,3) & & & & & \\ 
 \hline
30 & (0,3,16,21,23) & (2,20,10,0,25) & (0,27,4,10,11) & & mod $29$ & 174 & $\frac{3\times 29}{174}$\\
 & (1,0,15,2,26) & (9,0,4,12,26) & (20,8,$\infty$, 1,0) & & & &\\ 
 \hline
31 & (3,7,0,15,1) & (9,0,21,3,14) & (9,0,27,1,11) & & mod $31$ & 186 & $\frac{3\times 31}{186}$\\
 & (18,0,27,16,26) & (0,19,9,6,26) & (27,0,3,19,2) & & & &\\ 
  \hline
 35 & (6,7,0,30,12) & (0,24,18,10,31) & (0,4,6,5,21) & (23,0,32,8,25) & mod $34$ & 238 & $\frac{3\times 34+17}{238}$\\
 & (0,32,11,29,20) & (5,8,12,20,0) & (1,15,$\infty$, 4,0) & & & &\\ 
 \hline
36 & (1,0,6,9,21) & (3,13,0,2,27) & (13,22,0,18,26) & (2,17,9,30,0) & mod $36$  & 252 & $\frac{2\times 36+36+2\times 18}{252}$\\
& (4,5,34,0,16) & (10,14,0,11,30) &   & & & &\\ 
&(0,17,3,10,34)  & & & & & &\\
\hline

\end{tabular}
\end{center}
The above table has five columns. The first column contains the values of $v$ and the second column contains the base blocks. The third column shows that how to develope the elements of base blocks. Two last columns contain the number of blocks of corresponding design and the least possible value of $d$, respectively. For the remaining values of $v$, their associated super-simple directed designs are presented in the Appendix.
\end{proof}
\end{lemma}
\begin{lemma}\label{l2}
There exists a super-simple $(5,2)$-DGDD of type $5^5$ with $d\geq \frac{1}{2}$.
\begin{proof}
Let $X=Z_{25}$ and let $\mathcal{G}=\{\{i,5+i,10+i,15+i,20+i\}|  \  0\leq i\leq 4\}$. Here are the base blocks. These blocks are developed by +5 (mod 25).  
\begin{center}
\begin{tabular}{ccccc}
%\hline
%$v$ & base blocks & & & & $b_v$ & $d$ \\
 % \hline
 (4,0,22,21,23) & (1,18,22,24,0) & (11,2,9,0,23) & (2,0,8,6,24) & (9,7,18,0,6) \\
 (10,14,3,21,22) & (16,12,0,23,24) & (9,3,12,21,15) & (14,0,2,18,16) & (0,7,14,21,13) \\
 & & &  \\
 (1,23,14,2,10) & (4,6,10,13,12) & (0,7,1,4,3) & (0,9,13,1,12) & (3,6,12,14,0)  \\
 (3,6,20,17,9) &  (7,8,0,19,16) & (2,6,10,18,19) & (8,4,5,22,16) & (6,5,3,24,22) \\
 & &  &  & \\
%\hline
\end{tabular}
\end{center}
This directed group divisible design has $100$ blocks, contains $10$ disjoint directed trades of volume
2 in each of five columns. Since each
defining set for this design must contain one block of each directed trades, then each defining
set contains at least $50$ blocks. So $d\geq \frac{1}{2}$.
\end{proof}
\end{lemma}
\begin{lemma}\label{l3}
For each $t$, $6\leq t\leq 10$, there exists a super-simple $(5,2)$-DGDD of type $5^t$ with $d\geq \frac{1}{2}$.
\begin{proof}
Let the point set be $Z_{5t}$ and let the group set be $\{\{i,i+t,i+2t,i+3t,i+4t\}|  \ \  0\leq i\leq t-1\}$. The required base blocks are listed below. All the  base blocks are developed by mod $5t$.\\
\begin{center}
\begin{tabular}{|c|cccc|c|c|}
\hline
$t$ & base blocks & & & & $b_t$ & $d \geq$ \\
  \hline
 6 & (2,10,7,17,0) & (20,0,19,21,28) & (0,21,11,26,25) &  & $150$ & $\frac{2\times 30+30}{150}$\\
& (0,13,16,17,9) & (0,2,16,27,19) & & & & \\
  \hline
 7 & (0,1,31,33,30) & (19,0,10,25,27) & (0,16,4,22,33) &  & $210$ & $\frac{3\times 35}{210}$\\
& (11,0,15,20,23) & (9,0,19,22,18) & (0,1,20,12,25) & & & \\
\hline
 8 & (6,2,28,0,13) & (3,0,31,17,10) & (5,1,10,0,11) & (11,6,0,20,33) & $280$ & $\frac{3\times 40+20}{280}$ \\
& (2,23,0,4,1) & (22,0,5,25,28) & (0,12,27,2,31) & & & \\
\hline
9 & (16,17,0,30,32) & (0,3,20,41,4) & (11,0,26,37,32) & (2,0,33,30,37) & $360$ & $\frac{4\times 45}{360}$ \\
& (11,10,21,33,0) & (0,2,25,22,41) & (0,6,31,39,44) & (0,19,14,17,24) & & \\
\hline
10 & (12,21,0,33,25) & (12,3,11,0,26) & (0,13,17,36,28) & (0,3,31,47,49) & $450$ & $\frac{4\times 50+25}{450}$\\
& (18,2,23,45,0) & (16,0,11,17,35) & (2,5,38,0,46) & (21,28,0,2,37) & & \\
& & & & & & \\
& (0,6,43,7,32) & & & & & \\
\hline
\end{tabular}
\end{center}
\end{proof}
\end{lemma}
\begin{lemma}\label{l4}
There exists a super-simple $(5,2)$-DGDD of type $(15)^t$ for $t\in \{6,7,9\}$ with $d\geq \frac{1}{2}$.
\begin{proof}
Let the point set be $Z_{15t}$ and let the group set be $\{\{i,i+t,i+2t,\cdots,i+14t\}| \ \ 0\leq i\leq t-1\}$. The base blocks are listed below. Here, all the base blocks
are developed by mod $15t$ .
\begin{center}
\begin{tabular}{|c|cccc|c|c|}
\hline
$t$ & base blocks & & & & $b_t$ & $d \geq$ \\
  \hline
  6 & (32,15,0,83,1) & (16,56,15,0,43) & (85,2,88,81,0) & (0,77,1,57,82) & 1350 & $\frac{7\times 90+45}{1350}$\\
   & (0,22,68,85,33) & (40,87,37,59,0) & (26,51,55,64,0) & (63,55,0,76,23) & & \\
   & & & & & & \\
    & (44,0,41,7,15) & (0,31,47,51,70) & (0,45,73,71,44) &  (0,80,10,55,69) & & \\
    & (74,25,69,16,0) & (0,71,43,33,10) & (38,0,15,17,49) & & & \\
    \hline
7 & (0,31,57,75,76) & (0,50,55,102,72) & (27,40,93,0,99) & (0,66,74,82,93) & 1890 & $\frac{9\times 105}{1890}$\\
& (0,11,85,10,100) & (15,73,82,0,95) & (0,3,55,88,79) & (81,0,62,101,64) & & \\
& & & & & & \\
 &  (0,34,81,73,99) & (0,82,31,102,37) & (19,76,0,44,1) & (5,0,53,94,69) & & \\
& (76,22,0,58,68) & (0,51,92,54,94) & (45,0,43,81,93) & (0,69,92,96,101) & & \\
& & &  & &  & \\
 & (22,26,60,0,59) & & & & & \\
& (0,46,71,44,61) & & &  & & \\
\hline
9 & (25,55,60,0,85) & (12,49,52,80,0) & (0,4,2,10,48) & (22,7,0,51,110) & 3240 & $\frac{12\times 135}{3240}$\\
 & (23,0,105,65,76) & (38,6,0,107,22) & (8,4,20,0,96) & (23,0,84,71,58) & & \\
 & & & & & &  \\
 & (46,0,7,116,85) & (3,20,0,13,46) & (88,28,69,35,0) & (33,52,0,120,91) & & \\
 & (70,0,28,59,49) & (0,11,32,134,66) & (8,40,57,16,0) & (6,0,40,92,26) & & \\
 & & & & & & \\
 & (0,24,98,104,25) & (46,123,0,125,61) & (102,44,85,0,14) & (1,2,0,5,24) & & \\
 & (41,3,56,0,70) & (98,0,74,23,87) & (97,35,92,14,0) & (5,0,56,98,118) & & \\
 \hline
\end{tabular}
\end{center}
\end{proof}
\end{lemma}
\section{Main Theorem}
In this section we try to find super-simple $(v, 5, 2)$DDs for some admissible values of $v$ by recursive constructions
presented in Section 2 and using super-simple DGDDs obtained in Section 3. 
\begin{lemma}\label{l5}
There exists a super-simple $(v,5,2)$DD for each $v\in \{20i+\eta|  \ 5\leq i\leq 9, \eta=0,1\}$ with $d\geq \frac{1}{2}$.
\begin{proof}
Using  a super-simple $(5,2)$-DGDD of type $5^t$ for $5\leq t\leq 9$ with $d\geq \frac{1}{2}$ obtained in  Lemmas \ref{l2} and \ref{l3} and applying Construction \ref{1} with a TD$(5,4)$ as an input design comming from Lemma \ref{l11111}, we obtain a super-simple $(5,2)$-DGDD of type $(20)^t$ with $d\geq \frac{1}{2}$. On the other hand by Lemma \ref{l1}  there exists a super-simple $(20+\eta,5,2)$DD. So by Construction \ref{2} we obtain a super-simple $(20t+\eta,5,2)$DD with $d\geq \frac{1}{2}$, where $\eta= 0$ or $1$.
\end{proof}
\end{lemma}
\begin{lemma}\label{l6}
There exists a super-simple $(v,5,2)$DD for each $v\in \{125, 126, 145, 146, 150, 151\}$ with $d\geq \frac{1}{2}$.
\begin{proof}
We delete $5-a$ points from the last group of a TD$(6,5)$ coming from Lemma \ref{l11111} to obtain a $\{5,6\}$-GDD of type $5^5a^1$. Applying Construction \ref{1} and using a super-simple $(5, 2)$-DGDD of group type $5^5$ and $5^6$ with $d\geq \frac{1}{2}$ from Lemmas \ref{l2} and \ref{l3} we get  a super-simple $(5,2)$-DGDD of type $(25)^5(5a)^1$ with $d\geq \frac{1}{2}$. Since by Lemma \ref{l1} there exists a super-simple $(25+\eta,5,2)$DD
and a super-simple $(5a+\eta,5,2)$DD for $a\in \{0, 4, 5\}$ and $\eta=0,1$, by Construction \ref{2} we get a super-simple $(125+5a+\eta,5,2)$DD with $d\geq \frac{1}{2}
$.
\end{proof}
\end{lemma}
\begin{lemma}\label{l7}
There exists a super-simple $(v,5,2)$DD for $v\in \{155, 156\}$ with $d\geq \frac{1}{2}$.
\begin{proof}
Starting from a $5$-GDD of type $3^8 7^1$ (exists by Lemma $4.3$ in \cite{chen4}) and applying Construction \ref{1} by using a super-simple $(5,2)$-DGDD of type $5^5$ with $d\geq \frac{1}{2}$ as an input designs, we get a super-simple $(5,2)$-DGDD of type $(15)^8 (35)^1$ with $d\geq \frac{1}{2}$. Since by Lemma \ref{l1} there exists a super-simple $(15+\eta,5,2)$DD and a super-simple
$(35+\eta,5,2)$DD, by Construction \ref{2} we get a super-simple $(155+\eta,5,2)$DD with $d\geq \frac{1}{2}$, where $\eta= 0$ or $1$.
\end{proof}
\end{lemma}
\begin{lemma}\label{l8}
There exists a super-simple $(v,5,2)$DD for each $v\in \{170, 171, 175, 176, 185, 186\}$ with $d\geq \frac{1}{2}$.
\begin{proof}
Starting  from a $\{5,6,7,8\}$-GDD of type $4^76^1$ (exists by Lemma $4.4$ in \cite{chen4}) and  applying Construction \ref{1} by using a super-simple $(5,2)$-DGDD of type $5^5$, $5^6$, $5^7$ and $5^8$ with $d\geq \frac{1}{2}$ coming from Lemmas \ref{l2} and \ref{l3} we get a super-simple $(5,2)$-DGDD  of type $(20)^7 (30)^1$ with $d\geq \frac{1}{2}$. Since by Lemma \ref{l1} there exists a super-simple $(20+\eta,5,2)$DD and a super-simple $(30+\eta,5,2)$DD with $d\geq \frac{1}{2}$, then by Construction \ref{2} we get a super-simple $(170+\eta,5,2)$DD with $d\geq \frac{1}{2}$, where $\eta= 0$ or $1$.\\

Starting from a $TD(5,7)$ coming from Lemma \ref{l11111} and applying Construction \ref{1} by using a super-simple $(5,2)$-DGDD of type $5^5$ with $d\geq \frac{1}{2}$ coming from Lemma \ref{l2} we get a super-simple $(5,2)$-DGDD of type $(35)^5$ with $d\geq \frac{1}{2}$. Since by Lemma \ref{l1}
there exists a super-simple $(35+\eta,5,2)$DD, by Construction \ref{2} we get a super-simple $(175+\eta,5,2)$DD with $d\geq \frac{1}{2}$, where $\eta= 0$ or $1$.\\

Starting from a $\{5,6,7\}$-GDD of type $5^67^1$ (exists from Lemma $4.4$ in \cite{chen4}) and applying Construction \ref{1} by using a super-simple $(5,2)$-DGDD of type $5^5$, $5^6$ and $5^7$ with $d\geq \frac{1}{2}$ coming from Lemmas \ref{l2} and \ref{l3} we get a super-simple $(5,2)$-DGDD of type $(25)^6(35)^1$ with $d\geq \frac{1}{2}$. Since by Lemma \ref{l1} there exists a super-simple
$(25+\eta,5,2)$DD and a super-simple $(35+\eta,5,2)$DD, by Construction \ref{2} we get a super-simple $(185+\eta,5,2)$DD with $d\geq \frac{1}{2}$, where $\eta=0$ or $1$.
\end{proof}
\end{lemma}
\begin{lemma}\label{l9}
There exists a super-simple $(v,5,2)$DD for any $v\in \{90,91,105, 106, 135, 136\}$ with $d\geq \frac{1}{2}$.
\begin{proof}
By Lemma \ref{l4} there exists a super-simple $(5,2)$-DGDD of type $(15)^t$  with $d\geq \frac{1}{2}$ for $t\in\{6,7,9\}$. Since by Lemma \ref{l1} there exist a super-simple $(15+\eta,5,2)$DD with $d\geq \frac{1}{2}$ for $\eta=0,1$, by Construction \ref{2} we get a super-simple $(15t+\eta,5,2)$DD with $d\geq \frac{1}{2}$, where $\eta=0$ or $1$.
\end{proof}
\end{lemma}
\begin{lemma}\label{l10}
There exists a super-simple $(96,5,2)$DD with $d\geq \frac{1}{2}$.
\begin{proof}
A super-simple $(5, 2)$-DGDD of group type $4^6$ is listed as follows. Let $X = Z_{24}$ and 
$G = \{\{i, i+6,12+i,18+i\} | \ \
0\leq i\leq 5\}$. Below are the
required base blocks. All the  base blocks are developed by mod $24$. 
\begin{center}
\begin{tabular}{cc}
(0,2,1,4,11) & (1,0,5,22,15) \\
(13,2,0,16,21) & (0,1,20,9,16)\\
\end{tabular}
\end{center}
This super-simple DGDD has $96$ blocks, each of two columns has $24$ disjoint directed trades of volume $2$. Therefore each defining set for this super-
simple DGDD contains at least $24\times 2=48$ blocks. So $d\geq \frac{1}{2}$.\\
Starting from this DGDD and applying Construction \ref{1} with a $TD(5,4)$ coming from Lemma \ref{l11111} we get a super-simple $(5,2)$-DGDD of type $(16)^6$ with $d\geq \frac{1}{2}$. Since by Lemma \ref{l1} there exists a super-simple $(16,5,2)$DD with $d\geq \frac{1}{2}$, by Construction \ref{2} we get a
super-simple $(96,5,2)$DD with $d\geq \frac{1}{2}$.
\end{proof}
\end{lemma}
\begin{lemma}\label{l11}
There exists a super-simple $(v,5,2)$DD for any $v\in \{165,166\}$ with $d\geq \frac{1}{2}$.
\begin{proof}
Let the point set be $X = Z_{33}$ and the group set be
$G = \{\{i, i+11,i+22\} | \ \
0\leq i\leq 10\}$. Below are the
required base blocks. All the base blocks are developed by mod $33$.
\begin{center}
\begin{tabular}{ccc}
(6,2,0,3,27) & (10,0,26,2,19) & (1,0,4,6,5) \\
(9,15,19,0,29) & (1,13,20,0,8) & (2,0,15,30,5) \\
\end{tabular}
\end{center}
This super-simple DGDD has $198$ blocks, each of three columns has $33$ disjoint directed trades of volume $2$. Therefore each defining set for this super-
simple DGDD contains at least $33\times 3=99$ blocks. So $d\geq \frac{1}{2}$.\\
Starting from this DGDD and applying Construction \ref{1} with a $TD(5, 5)$ coming from Lemma \ref{l11111} we get  a super-simple $(5,2)$-DGDD of type $(15)^{11}$ with $d\geq \frac{1}{2}$. Since by Lemma \ref{l1} there exists a super-simple $(15+\eta,5,2)$DD, by Construction \ref{2} we obtain a super-simple $(165+\eta,5,2)$DD with $d\geq \frac{1}{2}$, where $\eta=0$ or $1$.
\end{proof}
\end{lemma}
\begin{lemma}\label{l12}
Suppose that $5\leq k\leq 10$ is an integer. Let $N(m)\geq k-2$, $r=k-5$ and let $M=\{5m,5a_1,\cdots,5a_r\}$, where $a_i\in [3,m]\cup \{0\}$, $1\leq i\leq r$. If there exists a super-simple $(l+\eta,5,2)$DD with $d\geq \frac{1}{2}$ for each $l\in M$, then there exists a super-simple $(25m+5\sum_{i=1}^r a_i+\eta,5,2)$DD with $d\geq \frac{1}{2}$, where $\eta=0$ or $1$. 
\begin{proof}
By Lemma $4.8$ in \cite{chen4}, there exists a
$\{5, 6, . . . , k\}$-GDD of type $m^5(a_1)^1(a_2)^1 \cdots  (a_r)^1$. Starting from this GDD and applying Construction \ref{1} by using a  super-simple $(5,2)$-DGDDs of type $5^t$ for $t\in \{5,6,\cdots, k\}$ with $d\geq \frac{1}{2}$ coming from Lemmas \ref{l2} and \ref{l3} we get a super-simple $(5,2)$-DGDD of type $(5m)^5(5a_1)^1(5a_2)^1\cdots (5a_r)^1$ with $d\geq \frac{1}{2}$. Since there exists a super-simple $(u+\eta,5,2)$DD for any
$u\in M$, by Construction \ref{2} we get a super-simple $(25m+5\sum_{i=1}^r a_i+\eta,5,2)$DD with $d\geq \frac{1}{2}$ where $\eta=0$ or $1$. 
\end{proof}
\end{lemma}
\begin{lemma}\label{l13}
There exists a super-simple $(v,5,2)$DD for any $v\in [190, 1591]^{0,1}_
5$ with $d\geq \frac{1}{2}$.
\begin{proof}
Applying Lemma \ref{l12} with parameters in the following table, we obtain a super-simple $(v, 5, 2)$DD for every
$v\in [190, 1591]^{0,1}_
5$ . All required $TD(k,m)$ exist by Lemma \ref{l11111}.
\begin{center}
\begin{tabular}{cccccc}
\hline
$v=25m+5\sum_{i=1}^r a_i+\eta$ \ \ \ \ \  & $m$ \ \ \ \ \ & $k$ \ \ \ \ &  $\sum_{i=1}^{k-5} a_i$ \ \ \ \  & $\eta$ \\
  \hline
 $ [190,281]^{0,1}_5$ \ \ \ \ \ \ & $7$ \ \ \ \ \ \ \ & $8$ \ \ \ & $[3,21]$ & $\{0,1\}$ \\
  $ [285,451]^{0,1}_5$ \ \ \ \ \ & $9$ \ \ \ \  & $10$ & \ \ \ $[12,45]$ & $\{0,1\}$ \\
   $ [455,651]^{0,1}_5$ \ \ \ \ \ & $13$ \ \ \ \ & $10$ & \ \ \ $[26,65]$ & $\{0,1\}$ \\
    $ [655,1251]^{0,1}_5$ \ \ \ \ \ & $25$ \ \ \ \ & $10$ & \ \ \ $[6,125]$ & $\{0,1\}$ \\
     $ [1255,1591]^{0,1}_5$ \ \ \ \ \ & $36$ \ \ \ \ & $10$ & \ \ \  $[71,138]$ & $\{0,1\}$ \\
     \hline
  \end{tabular}
\end{center}

\end{proof}
\end{lemma}
Now, we are in a position to conclude the main result.\\
\textit{{\bf Main Theorem.} For all $v\equiv 0,1$ (mod $5$) and $v\geq 15$, there exists a super-simple $(v,5,2)$DD with $d\geq \frac{1}{2}$. }
\begin{proof}
The proof is by induction on $v$. By  the above lemmas, the result is true for $v\in [15,1591]_{5}^{0,1}$.  Therefore, we assume that $v\geq 1595$. We can write $v=25m+5(a_1+a_2)+\eta$, where $m\geq 63$, $\eta=0,1$, $\{a_1,a_2\}\subset [3,m]\cup \{0\}$ and $a_1+a_2\in [3,2m]$. By induction
there exists a super-simple $(5m+\eta, 5, 2)$DD and a super-simple $(5a_i +\eta, 5, 2)$DD, for $i =1, 2$ with $d\geq \frac{1}{2}$. Since $N(m)\geq 5$, we know that there
exists a super-simple $(v, 5, 2)$DD by Lemma \ref{l12}.
\end{proof}
\section{Constructing of some LDPC codes by using super-simple $ (v,5,2) $DDs}

A low density parity check (LDPC) codes first proposed by Gallager in 1960 in his dissertation after that these codes were ignored about 36 years and then rediscovered by Mackey \cite{gallager,mackay}.
% LDPC codes with a carefully chosen parity check matrix have good distance
%properties \cite{gallager,mackay}.
\\
An $m\times n$ sparse matrix is one that many of its elements are equal to zero. An LDPC code is a linear code with a sparse parity check matrix.
\\
A binary code with parity check matrix with $d_c$ 1s in each column and $d_r$ 1s in
each row is a regular (or bi-regular) LDPC code. Otherwise it is an irregular LDPC
code.
\\
There is a natural association of an $m\times n$ binary matrix with a bipartite graph $ G = (L\cup R,E) $,
called Tanner graph \cite{tanner}, whose adjacency matrix is the parity check matrix of the code where columns are identified with the left nodes $ (L) $ of the bipartite
graph (the message nodes) and rows with the right nodes $ (R) $ of the graph (check nodes) and $ E $ is the set of edges.
%The variables associated with the check nodes are just the exclusive-OR of the variables
%associated with neighbor message nodes.
The girth of the code is the length of the shortest cycle in the Tanner graph $ G $ denoted by $ g $.
%The sparsity of the matrix leads to a graph where each node has relatively fewneighbors.
\\
Reducing the parity check matrix of an LDPC code to standard form to allow systematic coding may render it nonsparse, leading to inefficient encoding.
Better encoding techniques are given in \cite{rich}.
\\
Some approaches to construct LDPC codes are algebric-based, protograph-based, and convolutional LDPC codes.
Numerous techniques for the construction of LDPC codes have been proposed. These include the original codes of Gallager \cite{gallager}, MacKay codes \cite{mackay}, irregular
degree sequence codes, and codes based on combinatorial structures such as finite
geometries and designs \cite{amirzade, deng, gruner, vasic}.
\\
%Regular LDPC codes are used for the design of redundant arrays of independent disks, or RAID storage systems. Hellerstein et al. \cite{{helle}} observe that bi-regularity
%is desirable to minimize the cost of the basic read and write operations while minimizing loss by erasure. They employ Kirkman triple systems to construct such codes
%for RAID architectures permitting two erasures. This is generalized in \cite{chee} using
%t-designs and Tur´an systems.
%\\
%An edge incident with a message (check) node is of degree $i$ if it is incident with a message (check) node of degree $i$.
In this section, We use super-simple $ (v,5,2) $ directed designs to obtain parity-check matrices of trade-based LDPC codes, also we use a graph of trade in which its vertices are associated to the blocks of a super-simple directed design and each edge between two blocks shows a trade of volume 2, and then provide an example to show this method in an irregular LDPC code with girth 8.\\
%Let $ T=T_1-T_2 $ be a directed trade and  $\mathcal{D}=(V,\mathcal{B})$ be a super-simple directed design.
% If $ T_1 $ has some blocks of $ B $, $ T_2 $ obtains from $ T_1 $ by interchanging the order of common elements in blocks of $ T_1 $, it means
%we suppose that $ \mathcal{D} $ contains the directed trade $ T $.

Let $ V = \lbrace0,1,\ldots,v-1\rbrace $ be a $ v$-set, corresponding to a super-simple $ (v,5,2) $DD with $ n $ blocks $ b_1,b_2,\ldots,b_n $ we construct a $ (_{2}^{v}) \times n $ matrix $ A $ whose columns indices $ b_1,b_2,\ldots,b_n $ and rows indices $ (x_i,x_j) $, where $ x_i < x_j $ and $ x_i,x_j \in \lbrace0,1,\ldots,v-1\rbrace $, as follows:
\begin{center}
$A_{(x_i,x_j) l} = \begin {cases}

1 & if\quad (x_i,x_j)\in b_l $ which $ b_l $ appears in a trade,$ \\ 0 & o.w.

\end {cases}$
\end{center}
Then by removing all zero rows and columns of $A$, the resulting matrix denoted by $ C $ is considered as the parity-check matrix of trade-based LDPC code, if the number of rows of $C$ is less than or equal to the number of its columns. otherwise, $ C^\bot $ which is the transpose of $ C $, is taken as the parity-check matrix of the code. This parity-check matrix is the adjacency matrix of the Tanner graph.\\
The number of 1s in each rows and columns is at most $ 2\lambda $ and $ \lambda(k-1) $, respectively. Because any pair $ (x_i,x_j) $ where $ x_i < x_j $, occurs in $ \lambda $ blocks and each block of length $ k $ contains $ k-1 $ pairs $ (x_i,x_j) $s with $ x_i < x_j $, which may appear in $ \lambda $ trades.\\
Main theorem in previous section, shows that for all $v\equiv 0,1$ (mod $5$) and $v\geq 15$, there exists a super-simple $(v,5,2)$DD whose smallest defining sets have at least half
of the blocks. Because of this property we can construct parity-check matrices of LDPC codes in which
the number of rows and columns is equal to or bigger than the number of blocks. Thus, the length
of our proposed LDPC code is at least the number of blocks.\\
In the super-simple $ (v,5,2) $DD, since $ \lambda=2 $, three cases may occur for its directed trades:\\
$ 1. $ Two blocks have a common pair in which first block contains $ (x_i,x_j) $ and the second one contains $ (x_j,x_i) $, in this case the directed trade is of the form
\begin{center}
\begin{tabular}{cc|cc}
&$T_1$ & $T_2$ \\
\hline
&$b_1: (x_i,x_j, \ldots)$ & $( x_j,x_i,\ldots) $\\
&$b_2: ( x_j,x_i,\ldots) $ & $ (x_i,x_j, \ldots)$
\end{tabular}
\end{center}
$ 2. $ Three blocks have a common pair in which two blocks contain $ (x_i,x_j) $ and the third one contains $ (x_j,x_i) $, or vice versa. In this case we have two directed trades of volume 2 of the form\\
\begin{center}
\begin{tabular}{cc|cc}
&$T_1$ & $T_2$ \\
\hline
&$b_1: (x_i,x_j, \ldots)$ & $( x_j,x_i,\ldots) $\\
&$b_3: ( x_j,x_i,\ldots) $ & $ (x_i,x_j, \ldots)$
\end{tabular}
{\hspace{1cm}
\begin{tabular}{cc|cc}
&$T_1$ & $T_2$ \\
\hline
&$b_2: (x_i,x_j, \ldots)$ & $( x_j,x_i,\ldots) $\\
&$b_3: ( x_j,x_i,\ldots) $ & $ (x_i,x_j, \ldots)$
\end{tabular}
}
\end{center}
the third block that has $ (x_j,x_i) $ is common in two trades.\\
$ 3. $ Four blocks have a common pair in which two blocks contains $ (x_i,x_j) $ and two other ones contains $ (x_j,x_i) $, in this case we have four directed trades of volume 2 as shown below:\\
{\hfill
\begin{tabular}{cc|cc}
&$T_1$ & $T_2$ \\
\hline
&$b_1: (x_i,x_j, \ldots)$ & $( x_j,x_i,\ldots) $\\
&$b_3: ( x_j,x_i,\ldots) $ & $ (x_i,x_j, \ldots)$
\end{tabular}
}
{\hfill
\begin{tabular}{cc|cc}
&$T_1$ & $T_2$ \\
\hline
&$b_1: (x_i,x_j, \ldots)$ & $( x_j,x_i,\ldots) $\\
&$b_4: ( x_j,x_i,\ldots) $ & $ (x_i,x_j, \ldots)$
\end{tabular}
}
{\hfill
\begin{tabular}{cc|cc}
&$T_1$ & $T_2$ \\
\hline
&$b_2: (x_i,x_j, \ldots)$ & $( x_j,x_i,\ldots) $\\
&$b_3: ( x_j,x_i,\ldots) $ & $ (x_i,x_j, \ldots)$
\end{tabular}
}
{\hfill
\begin{tabular}{cc|cc}
&$T_1$ & $T_2$ \\
\hline
&$b_2: (x_i,x_j, \ldots)$ & $( x_j,x_i,\ldots) $\\
&$b_4: ( x_j,x_i,\ldots) $ & $ (x_i,x_j, \ldots)$
\end{tabular}
}\\
%Let $ V = \lbrace0,1, \cdots , v-1 \rbrace $, we construct binary matrix whose columns (message nodes) labeled with $ b_1 , \cdots b_n $ and rows (check nodes) labeled with pairs $ (x_i,x_j) $ where $ x_i < x_j \in V $, belongs to all blocks which appear in trades. This parity-check matrix defines as below:

%\begin{center}
%$C_{(x_i,x_j) l} = \begin {cases}

%1 & if \quad (x_i,x_j)\in b_l \\ 0 & o.w.

%\end {cases}$
%\end{center}

%If the number of columns of $ C $ is more than the number of rows then $ C $ is the parity-check of the trade-based LDPC code, otherwise $ C^{\perp} $ is the parity-check matrix of code. 
Next proposition describes the relation between the smallest volume of cyclical trade and the girth of trade-based LDPC code, which is used to construct the codes.
\begin{proposition} \label{p2}
\cite{amirzade} A super-simple directed design has a cyclical trade of volume s if and only if the
Tanner graph of the corresponding trade-based LDPC code has 2s-cycles.
\end{proposition}
\begin{proposition}\label{p3}
The girth of an LDPC code constructed by using super-simple $ (v,5,2) $DD is at least $ 6 $. 
\end{proposition}
\begin{proof}
From proposition \ref{p2}, we conclude that the girth of the Tanner graph of the corresponding  trade-based LDPC code constructed by using super-simple $ (v,5,2) $DD must be even. By the concept of super-simple directed designs every two blocks intersects in at most two points or an $ (x_i,x_j) $, therefore the Tanner graph is 4-cycle free, in other words there is no cycle in its Tanner graph of the following type:

\begin{figure}[h!]
\begin{center}
\begin{tikzpicture} [scale=.5, thin]
\vertex[fill] (B31) at (0,0) [label=above:\tiny $B_{i}$] {};
\vertex[fill] (B2) at (2,0) [label=above:\tiny $B_{j}$] {};
\vertex[fill] (B30) at (2,-2) [label=below:\tiny ] {};
\vertex[fill] (B14) at (0,-2) [label=below:\tiny] {};
\fill[] (2,-2)  circle (0mm) node[yshift=-.35cm]{\tiny $(x_i,x_j)$};
\fill[] (0,-2)  circle (0mm) node[yshift=-.35cm]{\tiny $(y_i,y_j)$};
\path
(B31) edge (B30)
(B2) edge (B30)
(B31) edge (B14)
(B2) edge (B14);
\end{tikzpicture}
\end{center}
\end{figure}

If we have cyclical trade of volume 3 results in 6-cycle in the Tanner graph of a trade-based LDPC code, otherwise the code has girth at least 8.
\end{proof}
Combining proposed technique to construct a parity-check matrix for an LDPC code based on the concept of trades in super-simple directed designs that explained above, proposition \ref{p2} and proposition \ref{p3}, we have the following result.
\begin{theorem}
The existence of a super-simple $ (v,5,2) $DD whose smallest defining set have at least half of the blocks can deduce an LDPC code with the girth at least 6.
\end{theorem} 

We define a special set that knowledge of this is important for a given LDPC code. We show the graphical structure of this set corresponding to the LDPC code in the next example.
\begin{definition}\cite{lin}
Let $ G $ be the Tanner graph of an LDPC code given by the null space of a parity-check matrix $ C_{m\times n} $. For $ 1\leq \alpha \leq n $ , $ 1\leq \beta \leq m $, an $ (\alpha , \beta) $-trapping set is a subset $ D\subseteq V $ of variable nodes in the Tanner graph of LDPC code such that $ \vert D\vert=\alpha $ and $\vert O(D)\vert= \beta $, where $ O(D) $ is the subset of check nodes of odd degree in an induced subgraph of the tanner graph $ G_D $. An elementary $ (\alpha,\beta) $-trapping sets ($ (\alpha,\beta) $ ETSs) are those trapping sets  for which all check nodes in the induced subgraph of the Tanner graph is of degree 1 or 2.
\end{definition}
In an $ (\alpha,\beta) $ ETS, by removing all check nodes of odd degree and replacing every check nodes of even degree with an edge we obtain a graph with $ \alpha $ vertices which is called normal graph.\\
%since every cycles in bipartite graph are even cycles, any cycles is presented by alternating sequence of check nodes and message nodes in Tanner graph. The length of the shortest cycle is called girth of the graph. 
There is a $2s$-cycle in the ETS if and only if there is a $s$-cycle in its corresponding normal graph. We use this fact to find easily the girth of the Tanner graph.\\
\begin{example}
The followings are the blocks of a super-simple $ (15,5,2) $DD in lemma \ref{l1} :
\begin{align*}
B_1&=(1, 0, 2, 3, 8)        \quad  &B_{15}&=(0, 1, 4, 10, 9)   \quad  &B_{29}&=(1, 0, \infty, 5, 7)  
\\                                                                    
B_2&=(3, 2, 4, 5, 10)       \quad  &B_{16}&=(3, 2, 6, 12, 11)  \quad  &B_{30}&=(2, 3, \infty, 7, 9)  
\\                                                                    
B_3&=(5, 4, 6, 7, 12)       \quad &B_{17}&=(5, 4, 8, 0, 13)    \quad  &B_{31}&=(4, 5, \infty, 9, 11) 
\\                                                                    
B_4&=(7, 6, 8, 9, 0)        \quad &B_{18}&=(7, 6, 10, 2, 1)    \quad  &B_{32}&=(6, 7, \infty, 11, 13)
\\                                                                    
B_5&=(9, 8, 10, 11, 2)      \quad &B_{19}&=(9, 8, 12, 4, 3)    \quad  &B_{33}&=(8, 9, \infty, 13, 1) 
\\                                                                    
B_6&=(11, 10, 12, 13, 4)    \quad &B_{20}&=(11, 10, 0, 6, 5)   \quad  &B_{34}&=(10, 11, \infty, 1, 3)
\\                                                                    
B_7&=(13, 12 ,0, 1, 6)      \quad &B_{21}&=(13, 12, 2, 8, 7)   \quad  &B_{35}&=(12, 13, \infty, 3, 5)    
\end{align*}
\begin{align*}
B_8&=(0, 3, 13, 11, 9)\quad &B_{22}&=(0, 7, 11, 4, 2) \quad  &B_{36}&=(13, 2, \infty, 0, 10) 
\\
B_9&=(2, 5, 1, 13, 11)\quad &B_{23}&=(2, 9, 13, 6, 4)\quad   &B_{37}&=(1, 4, \infty, 2, 12) 
\\ 
B_{10}&=(4, 7, 3, 1, 13) \quad &B_{24}&=(4, 11, 1, 8, 6) \quad   & B_{38}&=(3, 6, \infty, 4, 0)  
\\                          
B_{11}&=(6, 9, 5, 3, 1)\quad &B_{25}&=(6, 13, 3, 10, 8)\quad  &B_{39}&=(5, 8, \infty, 6, 2)   
\\                          
B_{12}&=(8, 11, 7, 5, 3) \quad &B_{26}&=(8, 1, 5, 12, 10) \quad  &B_{40}&=(7, 10, \infty, 8, 4)  
\\                             
B_{13}&=(10, 13, 9, 7, 5) \quad &B_{27}&=(10, 3, 7, 0, 12) \quad  &B_{41}&=(9, 12, \infty, 10, 6)
\\                             
B_{14}&=(12, 1, 11, 9, 7) \quad &B_{28}&=(12, 5, 9, 2, 0) \quad    &B_{42}&=(11, 0, \infty, 12, 8)
\end{align*}
As shown in Fig. \ref{f1}, each block appeares in $ \lbrace 2, 3, 4, 5\rbrace $ trades of volume 2 (We show each $ B_i $ by $ i$, $  i \in \lbrace1,2,\cdots42\rbrace $ in Fig. \ref{f1}). The graph is free of triangles which showes that there is no cyclical trades of volume 3. In the normal graph the smallest cycle is of size 4 which corresponds to a $CT_4$. By using this design we can construct a matrix $C$ which has the same $42$ rows and columns. Therefore $C$ is the parity-check matrix of an irregular LDPC code with 
%column weights in 
$d_c=\lbrace 1, 2, 3, 4 \rbrace$ ones in columns and %row weights in 
$d_r=\lbrace 2, 3, 4 \rbrace$ ones in rows. Fig. 2 corresponds to the tanner graph of irregular LDPC code of super-simple  $(15,5,2)$ DD which the girth of it is 8, this girth that shown in Fig. 3 can be $g =\lbrace B_2,(4,5),B_{31},(9,11),B_{14},(7,9),B_{30},(2,3),B_2 \rbrace $.\\
\begin{figure}
\centering
\includegraphics[scale=0.6]{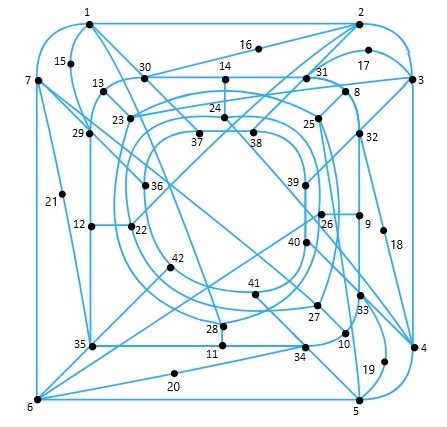}
\caption{The graph corresponding to directed trades of super-simple $(15,5,2)$ DD}
\end{figure}\label{f1}
The same can be done for other existing super-simple $ (v,5,2) $DD.
\begin{figure}[h!]
\begin{flushleft}
\begin{tikzpicture} [scale=.5, thin]
\vertex[fill] (B1) at (10,0) [label=above:\tiny $B_{1}$] {};
\vertex[fill] (B2) at (11,0) [label=above:\tiny $B_{2}$] {};
\vertex[fill] (B3) at (12,0) [label=above:\tiny $B_{3}$] {};
\vertex[fill] (B4) at (13,0) [label=above:\tiny $B_{4}$] {};
\vertex[fill] (B5) at (14,0) [label=above:\tiny $B_{5}$] {};
\vertex[fill] (B6) at (15,0) [label=above:\tiny $B_{6}$] {};
\vertex[fill] (B7) at (16,0) [label=above:\tiny $B_{7}$] {};
\vertex[fill] (B8) at (17,0) [label=above:\tiny $B_{8}$] {};
\vertex[fill] (B9) at (18,0) [label=above:\tiny $B_{9}$] {};
\vertex[fill] (B10) at(19,0) [label=above:\tiny $B_{10}$] {};
\vertex[fill] (B11) at(20,0) [label=above:\tiny $B_{11}$] {};
\vertex[fill] (B12) at(21,0) [label=above:\tiny $B_{12}$] {};
\vertex[fill] (B13) at(22,0) [label=above:\tiny $B_{13}$] {};
\vertex[fill] (B14) at(23,0) [label=above:\tiny $B_{14}$] {};
\vertex[fill] (B15) at(24,0) [label=above:\tiny $B_{15}$] {};
\vertex[fill] (B16) at(25,0) [label=above:\tiny $B_{16}$] {};
\vertex[fill] (B17) at(26,0) [label=above:\tiny $B_{17}$] {};
\vertex[fill] (B18) at(27,0) [label=above:\tiny $B_{18}$] {};
\vertex[fill] (B19) at(28,0) [label=above:\tiny $B_{19}$] {};
\vertex[fill] (B20) at(29,0) [label=above:\tiny $B_{20}$] {};
\vertex[fill] (B21) at(30,0) [label=above:\tiny $B_{21}$] {};
\vertex[fill] (B22) at (30,-1) [label=right:\tiny $B_{22}$] {};
\vertex[fill] (B23) at (30,-2) [label=right:\tiny $B_{23}$] {};
\vertex[fill] (B24) at (30,-3) [label=right:\tiny $B_{24}$] {};
\vertex[fill] (B25) at (30,-4) [label=right:\tiny $B_{25}$] {};
\vertex[fill] (B26) at (30,-5) [label=right:\tiny $B_{26}$] {};
\vertex[fill] (B27) at (30,-6) [label=right:\tiny $B_{27}$] {};
\vertex[fill] (B28) at (30,-7) [label=right:\tiny $B_{28}$] {};
\vertex[fill] (B29) at (30,-8) [label=right:\tiny $B_{29}$] {};
\vertex[fill] (B30) at (30,-9) [label=right:\tiny $B_{30}$] {};
\vertex[fill] (B31) at (30,-10) [label=right:\tiny $B_{31}$] {};
\vertex[fill] (B32) at (30,-11) [label=right:\tiny $B_{32}$] {};
\vertex[fill] (B33) at (30,-12) [label=right:\tiny $B_{33}$] {};
\vertex[fill] (B34) at (30,-13) [label=right:\tiny $B_{34}$] {};
\vertex[fill] (B35) at (30,-14) [label=right:\tiny $B_{35}$] {};
\vertex[fill] (B36) at (30,-15) [label=right:\tiny $B_{36}$] {};
\vertex[fill] (B37) at (30,-16) [label=right:\tiny $B_{37}$] {};
\vertex[fill] (B38) at (30,-17) [label=right:\tiny $B_{38}$] {};
\vertex[fill] (B39) at (30,-18) [label=right:\tiny $B_{39}$] {};
\vertex[fill] (B40) at (30,-19) [label=right:\tiny $B_{40}$] {};
\vertex[fill] (B41) at (30,-20) [label=right:\tiny $B_{41}$] {};
\vertex[fill] (B42) at (30,-21) [label=right:\tiny $B_{42}$] {};

\vertex[fill] (01) at (9,-1)  [label=left:\tiny {0,1}] {};
\vertex[fill] (02) at (9,-2)  [label=left:\tiny {0,2}] {};
\vertex[fill] (07) at (9,-3) [label=left:\tiny {0,7}] {};
\vertex[fill] (012) at(9,-4)  [label=left:\tiny {0,12}] {};
\vertex[fill] (0*) at (9,-5) [label=left:\tiny {$0$,$\infty$}] {};
\vertex[fill] (13) at (9,-6)  [label=left:\tiny {1,3}] {};
\vertex[fill] (15) at (9,-7) [label=left:\tiny {1,5}] {};
\vertex[fill] (18) at (9,-8)  [label=left:\tiny {1,8}] {};
\vertex[fill] (111) at(9,-9)  [label=left:\tiny {1,11}] {};
\vertex[fill] (113) at(9,-10) [label=left:\tiny {1,13}] {};
\vertex[fill] (23) at (9,-11)[label=left:\tiny {2,3}] {};
\vertex[fill] (24) at (9,-12) [label=left:\tiny {2,4}] {};
\vertex[fill] (29) at (9,-13)[label=left:\tiny {2,9}] {};
\vertex[fill] (2*) at (9,-14) [label=left:\tiny {$2$,$\infty$}] {};
\vertex[fill] (35) at (9,-15)[label=left:\tiny {3,5}] {};
\vertex[fill] (37) at (9,-16) [label=left:\tiny {3,7}] {};
\vertex[fill] (310) at(9,-17) [label=left:\tiny {3,10}] {};
\vertex[fill] (313) at(9,-18) [label=left:\tiny {3,13}] {};
\vertex[fill] (45) at (9,-19)[label=left:\tiny {4,5}] {};
\vertex[fill] (46) at (9,-20) [label=left:\tiny {4,6}] {};
\vertex[fill] (411) at(9,-21) [label=left:\tiny {4,11}] {};

\vertex [color=white] (333) at(3,-22)  {};

\vertex[fill] (4*) at (9,-22) [label=below:\tiny {$4$,$\infty$}] {};
\vertex[fill] (57) at (10,-22)  [label=below:\tiny {5,7}] {};
\vertex[fill] (59) at (11,-22)  [label=below:\tiny {5,9}] {};
\vertex[fill] (512) at(12,-22)  [label=below:\tiny {5,12}] {};
\vertex[fill] (67) at (13,-22)  [label=below:\tiny {6,7}] {};
\vertex[fill] (68) at (14,-22)  [label=below:\tiny {6,8}] {};
\vertex[fill] (613) at(15,-22)  [label=below:\tiny {6,13}] {};
\vertex[fill] (6*) at (16,-22)  [label=below:\tiny {$6$,$\infty$}] {};
\vertex[fill] (79) at (17,-22)  [label=below:\tiny {7,9}] {};
\vertex[fill] (711) at(18,-22)  [label=below:\tiny {7,11}] {};
\vertex[fill] (89) at (19,-22) [label=below:\tiny {8,9}] {};
\vertex[fill] (810) at(20,-22) [label=below:\tiny {8,10}] {};
\vertex[fill] (8*) at (21,-22) [label=below:\tiny{ $8$,$\infty$}] {};
\vertex[fill] (911) at(22,-22) [label=below:\tiny {9,11}] {};
\vertex[fill] (913) at(23,-22) [label=below:\tiny {9,13 }] {};
\vertex[fill] (1011)at(24,-22) [label=below:\tiny {10,11 }] {};
\vertex[fill] (1012)at(25,-22) [label=below:\tiny {10,12 }] {};
\vertex[fill] (10*) at(26,-22) [label=below:\tiny {10,$\infty$}] {};
\vertex[fill] (1113)at(27,-22) [label=below:\tiny {11,13}] {} ;
\vertex[fill] (1213)at(28,-22) [label=below:\tiny {12,13 }] {};
\vertex[fill] (12*) at(29,-22) [label=below:\tiny { $12$,$\infty$}] {};
\path
(01) edge (B1)
(01) edge (B7)
(01) edge (B15)
(01) edge (B29)

(02) edge (B1)
(02) edge (B28)

(07) edge (B22)
(07) edge (B27)

(012) edge (B7)
(012) edge (B27)

(0*) edge (B29)
(0*) edge (B36)
(0*) edge (B42)

(13) edge (B10)
(13) edge (B11)
(13) edge (B34)

(15) edge (B9)
(15) edge (B26)

(18) edge (B24)
(18) edge (B26)

(111) edge (B14)
(111) edge (B24)

(113) edge (B9)
(113) edge (B10)
(113) edge (B33)

(23) edge (B1)
(23) edge (B2)
(23) edge (B16)
(23) edge (B30)

(24) edge (B2)
(24) edge (B22)

(29) edge (B23)
(29) edge (B28)

(2*) edge (B30)
(2*) edge (B36)
(2*) edge (B37)

(35) edge (B11)
(35) edge (B12)
(35) edge (B35)

(37) edge (B10)
(37) edge (B27)

(310) edge (B25)
(310) edge (B27)

(313) edge (B8)
(313) edge (B25)

(45) edge (B2)
(45) edge (B3)
(45) edge (B17)
(45) edge (B31)

(46) edge (B3)
(46) edge (B23)

(411) edge (B22)
(411) edge (B24)

(4*) edge (B31)
(4*) edge (B37)
(4*) edge (B38)

(57) edge (B12)
(57) edge (B13)
(57) edge (B29)

(59) edge (B11)
(59) edge (B28)

(512) edge (B26)
(512) edge (B28)

(67) edge (B3)
(67) edge (B4)
(67) edge (B18)
(67) edge (B32)

(68) edge (B4)
(68) edge (B24)

(613) edge (B23)
(613) edge (B25)

(6*) edge (B32)
(6*) edge (B38)
(6*) edge (B39)

(79) edge (B13)
(79) edge (B14)
(79) edge (B30)

(711) edge (B12)
(711) edge (B22)

(89) edge (B4)
(89) edge (B5)
(89) edge (B19)
(89) edge (B33)

(810) edge (B5)
(810) edge (B25)

(8*) edge (B33)
(8*) edge (B39)
(8*) edge (B40)

(911) edge (B8)
(911) edge (B14)
(911) edge (B31)

(913) edge (B13)
(913) edge (B23)

(1011) edge (B5)
(1011) edge (B6)
(1011) edge (B20)
(1011) edge (B34)

(1012) edge (B6)
(1012) edge (B26)

(10*) edge (B34)
(10*) edge (B40)
(10*) edge (B41)

(1113) edge (B8)
(1113) edge (B9)
(1113) edge (B32)

(1213) edge (B6)
(1213) edge (B7)
(1213) edge (B21)
(1213) edge (B35)

(12*) edge (B35)
(12*) edge (B41)
(12*) edge (B42);
\end{tikzpicture}
\end{flushleft}\caption{The Tanner graph of irregular LDPC code of girth 8 corresponding to super-simple $ (15,5,2) $DD   }
\end{figure}\label{f2}
\begin{figure}[h!]
\begin{center}
\begin{tikzpicture} [scale=.5, thin]
\vertex[fill] (B31) at (0,0) [label=above:\tiny $B_{31}$] {};
\vertex[fill] (B2) at (2,0) [label=above:\tiny $B_{2}$] {};
\vertex[fill] (B30) at (2,-2) [label=below:\tiny $B_{30}$] {};
\vertex[fill] (B14) at (0,-2) [label=below:\tiny $B_{14}$] {};
\path
(B31) edge (B2)
(B2) edge (B30)
(B30) edge (B14)
(B14) edge (B31);
\hspace{5cm}
\vertex[fill] (B2) at (0,0) [label=above:\tiny $B_{31}$] {};
\vertex[fill] (B31) at (2,0) [label=above:\tiny $B_{2}$] {};
\vertex[fill] (B14) at (4,0) [label=above:\tiny $B_{30}$] {};
\vertex[fill] (B30) at (6,0) [label=above:\tiny $B_{14}$] {};
\vertex[fill] (45) at (1,-2) [label=below:\tiny ${4,5}$] {};
\vertex[fill] (911) at (3,-2) [label=below:\tiny ${9,11}$] {};
\vertex[fill] (79) at (5,-2) [label=below:\tiny ${7,9}$] {};
\vertex[fill] (23) at (7,-2) [label=below:\tiny ${2,3}$] {};
\vertex[fill] (B8) at (4,-3) [label=below:\tiny $B_{8}$] {};
\vertex[fill] (B13) at (6,-3) [label=below:\tiny $B_{13}$] {};
\vertex [color=white] (333) at(17,-3)  {};
\path
(B2) edge (45)
(B2) edge (23)
(B31) edge (911)
(B31) edge (45)
(B14) edge (79)
(B14) edge (911)
(B30) edge (23)
(B30) edge (79)
(911) edge (B8)
(79) edge (B13);
\end{tikzpicture}
\end{center}\caption{A normal graph ($ 4 $-cycle) and its corresponding $ (4,2) $ ETS ($ 8 $-cycle) }
\end{figure}\label{f3}
\end{example}
\newpage
\mbox{}
\section*{Appendix}
\begin{center}
\begin{tabular}{|c|cccc|c|c|c|}
\hline
$v$ & base blocks & & & & & $b_v$  & $d$ \\
  \hline
40 & (7,8,18,5,0) & (0,11,13,27,36) & (0,14,17,2,22) & (0,20,24,30,6) & mod 39 & 312 & $\frac{4\times 39}{312}$\\
& (0,38,22,12,30) & (0,28,$\infty$, 4,33) & (6,0,1,38,18)  &  (0,7,26,35,3) & & &\\ 
 \hline
41 & (0,4,1,11,29) & (6,8,27,0,32) & (0,11,7,10,23) & (0,39,20,6,15) & mod 41 & 328 & $\frac{4\times 41}{328}$\\
& (36,39,28,7,0) & (1,19,0,25,15) & (0,8,38,36,29)  &  (0,23,1,27,17) & & &\\
 \hline
45 & (0,11,21,2,7) & (0,40,9,26,43) & (0,18,38,32,39) & (0,3,1,29,19)  & mod 44 & 396 & $\frac{4\times 44+22}{396}$\\
& (0,4,10,23,12) & (0,41,24,5,36) & (3,9,$\infty$, 0,32)  &  (0,16,43,14,36) & & & \\ 
 & &  & & & & & \\
 & (22,0,9,33,37) & & & & & &\\
  \hline
46 & (14,1,7,0,10) & (12,0,24,26,32) & (31,0,28,25,29) & (0,5,24,22,45) & mod 46 & 414 & $\frac{4\times 46+23}{414}$\\
 & (0,7,2,16,33) & (10,0,15,37,28) & (17,0,25,35,38) & (0,20,1,36,31) & & & \\
 & & & & & & &\\
 & (23,0,30,42,34) & & & & & & \\ 
 \hline
50 & (0,15,8,47,48) & (5,45,0,19,42) & (1,23,0,3,31) & (0,5,25,12,36) & mod 49 & 490 & $\frac{5\times 49}{490}$\\
 & (0,32,41,11,45) & (0,30,45,47,16) & (0,38,35,41,10) & (20,0,27,33,21) & & &\\ 
 & & & & & & &\\
 & (5,0,17,27,43) & & & & & & \\
 & (0,5,$\infty$,14,39) & & & & & &\\
 \hline
51 & (40,17,0,22,29) & (17,27,30,32,0) & (0,33,20,39,50) & (5,27,35,0,42) & mod 51 & 510 & $\frac{5\times 51}{510}$\\
 & (15,0,18,47,31) & (6,42,0,32,44) & (0,4,31,41,8) & (28,2,42,16,0) & & &\\ 
 & & & & & & & \\
 & (24,30,0,50,47) & & & & & & \\
 & (7,8,5,50,0) & & & & & & \\
  \hline
 55 & (15,25,53,3,0) & (6,8,0,28,31) & (14,7,22,0,40) & (6,0,21,42,44) & mod 54 & 594 & $\frac{5\times 54+27}{594}$\\
 & (2,20,3,47,0) & (0,9,5,13,16) & (0,20,31,45,30) & (5,11,$\infty$,0,24) & & &\\ 
 & & & & & & & \\
 & (4,13,18,30,0) & (17,0,6,33,52) & & & & & \\
 & (5,0,35,22,34) & & & & & & \\
 \hline
56 & (4,29,0,10,48) & (0,40,21,39,41) & (0,23,15,20,28) & (0,7,29,6,46) & mod 56  & 616 & $\frac{2\times 56+4\times 56}{616}$\\
 & (0,9,27,13,30) & (4,2,32,0,36) & (0,51,11,44,26) & (6,20,13,32,0) & & &\\ 
 & (9,11,25,0,34) & & & & & & \\
 &  & (0,10,42,45,53) & & & & & \\
 &  & (0,46,12,47,41) & & & & & \\
 \hline
60 & (4,1,0,10,52) & (2,37,5,0,41) & (6,31,0,19,52) & (0,44,$\infty$,14,23) & mod 59 & 708 & $\frac{6\times 59}{708}$ \\
 & (3,0,1,17,34) & (0,3,8,21,32) & (0,5,16,42,6) & (16,0,4,40,31) & & &\\ 
 & & & & & & & \\
 & (8,0,2,20,45) & (25,10,0,32,12) & & & & & \\
 & (0,41,10,40,48) & (5,24,0,20,50) & & & & & \\
 \hline
 61 & (3,0,55,1,21) & (0,2,6,49,42) & (4,0,12,23,37) & (24,8,0,46,13) & mod 61 & 732  & $\frac{6\times 61}{732}$\\
 & (3,0,9,2,43) & (6,0,18,4,25) & (12,0,36,8,50) & (0,16,48,31,26) & & & \\ 
 & & & & & & & \\
 & (32,0,35,1,52) & (0,11,24,16,39) & & & & & \\
 & (35,0,3,44,34) & (0,48,22,32,17) & & & & & \\
 \hline
 65 & (10,25,0,23,39) & (23,20,0,22,48) & (0,23,19,53,51) & (0,27,36,35,5) & mod 64 & 832 & $\frac{6\times 64+32}{832}$\\
 & (0,41,$\infty$,31,52) & (0,10,22,27,1) & (0,46,25,19,49) & (0,37,51,31,57) & & & \\ 
 & & & & & & & \\
 & (0,7,45,16,56) & (4,49,44,0,56) & (28,14,0,32,61) & & & & \\
 & (5,47,7,0,53) & (9,56,0,10,13) & & & & & \\
 \hline
\end{tabular}
\end{center}
\begin{center}
\begin{tabular}{|c|cccc|c|c|c|}
\hline
$v$ & base blocks & & & & & $b_v$ &  $d$ \\
  \hline
66 & (17,37,36,0,65) & (15,27,5,0,24) & (22,4,0,20,27) & (10,0,17,25,43)  & mod 66 & 858 & $\frac{6\times 66+33}{858}$\\
& (20,0,31,62,24) & (6,29,63,0,42) & (5,0,35,57,9)  &  (3,13,0,49,54) & & &\\ 
& & & & & & & \\
& (34,0,6,32,59) & (45,0,58,39,47) & (32,16,0,33,54) & & & &\\
& (0,3,14,15,55) & (29,0,43,35,45) & & & & &\\
 \hline
70 & (0,26,33,52,64) & (5,59,0,25,43) & (18,37,47,0,65) & (30,35,0,37,36) & mod 69 & 966 & $\frac{7\times 69}{966}$\\
& (9,23,38,0,62) & (32,0,44,23,57) & (29,33,0,9,50)  &  (2,37,$\infty$,13,0) & & &\\
& & & & & & &\\
& (0,9,11,59,51) & (22,52,0,55,49) & (0,16,13,21,43) & & & &\\
& (4,27,0,20,28) & (6,12,0,4,58) & (1,0,15,56,59) & & & &\\
 \hline
 71 & (8,48,41,35,0) & (20,0,67,49,52) & (0,66,41,1,58) & (0,40,7,13,48) & mod 71 & 994 & $\frac{7\times 71}{994}$\\
& (0,7,4,27,42) & (47,0,69,19,37) & (1,6,0,31,14)  &  (23,27,0,20,56) & & &\\
& & & & & & &\\
& (43,0,55,34,45) & (34,0,62,50,60) & (19,0,43,53,21) & & & &\\
& (62,0,16,17,5) & (16,25,0,70,11) & (0,67,15,47,18) & & & &\\
 \hline
75 & (0,30,19,25,47) & (5,1,0,73,57) & (18,36,0,20,47) & (40,55,0,8,60) & mod 74 & 1110 & $\frac{7\times 74+37}{1110}$\\
& (0,44,36,51,57) & (37,0,4,33,55) & (0,36,21,60,61)  &  (0,23,31,16,28) & & &\\
& & & & & & &\\
& (0,30,39,10,65) & (24,0,72,31,65) & (0,32,1,49,46) & (12,25,28,65,0) & & &\\
& (32,21,64,0,44) & (0,4,$\infty$,66,68) & (0,3,14,38,64) & & & &\\
  \hline
76 & (36,17,58,0,71) & (5,58,7,0,65) & (0,5,51,55,67) & (15,0,27,42,56) & mod 76 & 1140 & $\frac{7\times 76+38}{1140}$\\
& (0,36,53,20,55) & (0,38,63,37,48) & (6,12,0,5,45)  &  (27,0,8,37,51) & & &\\
& & & & & & &\\
& (0,9,1,67,22) & (0,3,4,52,23) & (10,0,54,6,52) & (42,0,74,30,68) & & &\\
& (29,37,52,0,73) & (11,45,0,73,56) & (0,30,17,26,33) & & & &\\
  \hline
  80 & (21,38,20,0,28) & (32,8,33,0,11) & (0,50,44,51,40) & (0,49,64,4,40) & mod 79 & 1264 & $\frac{8\times 79}{1264}$\\
& (37,70,0,20,68) & (40,0,61,39,63) & (73,0,22,54,20)  &  (19,0,43,65,13) & & &\\
& & & & & & &\\
& (0,38,$\infty$,64,10) & (16,4,7,0,52) & (0,2,44,17,76) & (65,23,29,0,33) & & &\\
& (0,66,41,53,71) & (48,45,61,56,0) & (0,37,53,67,72) & (0,35,62,12,21) & & &\\
  \hline
  81 & (49,0,47,26,30) & (27,38,0,65,79) & (7,31,66,68,0) & (0,16,58,21,52) & mod 81 & 1296 & $\frac{8\times 81}{1296}$\\
& (0,41,15,49,31) & (0,50,7,75,78) & (8,0,7,64,18)  &  (37,42,0,74,46) & & &\\
& & & & & & &\\
& (52,0,17,60,78) & (14,0,20,68,69) & (0,33,45,67,73) & (0,56,40,69,70) & & &\\
& (9,18,0,45,57) & (28,48,0,51,70) & (57,0,59,80,76) & (9,34,19,4,0) & & &\\
  \hline
  85 & (4,7,60,0,82) & (48,0,68,25,11) & (25,41,0,57,58) & (0,18,68,30,47) & mod 84 & 1428 & $\frac{8\times 84+42}{1428}$\\
& (0,10,45,15,75) & (17,35,0,48,73) & (15,0,41,67,74)  &  (56,38,0,53,42) & & &\\
& & & & & & &\\
& (18,52,57,0,55) & (1,22,0,72,65) & (4,0,78,28,48) & (29,0,60,38,40) & & &\\
& (5,0,$\infty$,19,40) & (75,0,76,4,83) & (0,49,72,51,57) & (3,26,48,40,0) & & &\\
& & & & & & &\\
& (0,29,71,39,6) & & & & & &\\
  \hline
\end{tabular}
\end{center}
\begin{center}
\begin{tabular}{|c|cccc|c|c|c|}
\hline
$v$ & base blocks & & & & & $b_v$ & $d$ \\
  \hline
  86 & (68,41,72,69,0) & (0,59,84,31,33) & (0,63,43,7,78) & (0,22,77,21,79) & mod 86 & 1462  & $\frac{8\times 86+43}{1462}$\\
& (13,64,0,12,67) & (10,0,72,70,82) & (0,23,41,50,47)  &  (0,27,5,37,46) & & &\\
& & & & & & &\\
& (30,78,24,0,44) & (0,36,26,78,65) & (79,0,38,74,49) & (17,46,70,0,28) & & &\\
& (16,0,42,65,19) & (0,8,69,82,75) & (0,39,1,17,54) & (0,20,77,66,71) & & &\\
& & & & & & &\\
& (5,18,39,0,43) & & & & & &\\
  \hline
95 & (3,77,0,1,36) & (40,1,0,13,4) & (0,80,82,64,69) & (48,0,92,67,59)  & mod 94 & 1598 & $\frac{8\times 94+47}{1598}$\\
& (49,53,0,77,84) & (21,73,0,81,38) & (79,0,89,14,10)  &  (46,7,0,48,75) & & &\\ 
& & & & & & &\\
& (3,68,86,52,0) & (12,28,62,0,13) & (19,58,0,63,25) & (38,72,1,8,0) & & &\\
& (0,6,$\infty$,20,72) & (55,0,26,88,37) & (0,31,88,16,40) & (23,0,3,50,43) & & &\\
& & & & & & &\\
& (32,0,28,23,70) & & & & & &\\
 \hline
110 & (0,101,30,34,43) & (36,0,103,19,83) & (46,0,74,57,98) & (0,1,44,50,67) & mod 109 & 2398 & $\frac{11\times 109}{2398}$\\
& (0,18,26,40,98) & (40,85,0,88,97) & (40,38,0,31,48)  &  (4,70,81,14,0) & & &\\
& & & & & & &\\
& (0,100,73,78,85) & (0,33,87,46,53) & (0,15,108,94,75) & (0,59,77,79,81) & & &\\
& (23,48,63,0,104) & (65,8,0,55,84) & (3,19,$\infty$,54,0) & (86,37,0,107,102) & & &\\
& & & & & & &\\
& (0,74,71,108,61) & (21,0,86,59,89) & (51,26,0,96,80) & & & &\\
& (0,36,50,85,32) & (0,6,33,97,27) & (0,32,1,37,63) & & & &\\
 \hline
111 & (8,0,80,85,104) & (33,69,68,0,76) & (0,55,87,97,57) & (0,33,34,108,86) & mod 111 & 2442 & $\frac{11\times 111}{2442}$\\
& (2,19,50,41,0) & (0,68,5,105,79) & (57,47,0,88,103)  &  (44,28,104,0,110) & & &\\
& & & & & & &\\
& (0,3,20,107,86) & (51,0,63,91,64) & (23,0,52,96,39) & (79,6,27,77,0) & & &\\
& (53,36,56,0,80) & (5,0,51,107,41) & (0,59,15,40,77) & (38,83,0,61,65) & & &\\
& & & & & & &\\
& (11,0,93,79,33) & (0,69,81,95,99) & (33,20,82,0,90) & & & &\\
& (0,47,49,58,72) & (17,12,0,38,47) & (0,53,50,63,69) & & & &\\
\hline
115 & (0,106,11,1,109) & (0,113,53,111,102) & (46,16,0,74,1) & (50,43,30,0,42) & mod 114 & 2622 & $\frac{11\times 114+57}{2622}$\\
& (18,50,106,87,0) & (0,4,6,83,29) & (80,98,57,112,0)  &  (93,100,34,79,0) & & &\\
& & & & & & & \\
& (0,102,12,79,82) & (0,85,11,74,7) & (75,97,68,101,0) & (42,0,62,108,15) & & &\\
& (51,101,34,0,30) & (88,10,26,0,51) & (54,0,71,18,5) & (36,0,76,44,86) & & &\\
& & & & & & &\\
& (0,100,24,51,29) & (0,75,$\infty$,31,6) & (65,41,0,10,109) & (2,0,92,35,89) & & &\\
& (92,0,9,37,75) & (84,75,0,56,23) & (0,52,72,19,34) & & & &\\
\hline
116 & (1,94,67,88,0) & (0,101,63,54,68) & (11,16,0,35,67) & (63,0,106,39,115) & mod 116 & 2668 & $\frac{(6+2+4)\times 116}{2668}$\\
& (8,0,42,65,11) & (41,56,0,87,113) & (99,5,0,16,12)  &  (0,35,41,20,60) & & &\\
& & & & & & &\\
& (0,62,55,99,100) & (26,0,41,114,107) & (0,37,61,25,23) & (0,4,34,102,74) & & &\\
& (0,64,71,18,74) & (114,0,50,82,108) & (24,0,89,99,103) & (34,0,12,98,9) & & &\\
& &  & (66,68,95,74,0) & (0,47,77,90,85) & & &\\
& &  &  & (19,0,55,86,73) & & &\\
& &  &  & (12,0,106,51,71) & & &\\
& &  &  & (40,20,99,36,0) & & &\\
& &  &  & (33,69,0,1,46) & & &\\
& &  &  & (39,0,28,97,72) & & &\\
  \hline
\end{tabular}
\end{center}
\begin{center}
\begin{tabular}{|c|cccc|c|c|c|}
\hline
$v$ & base blocks & & & & & $b_v$ & $d$ \\
  \hline
  130 & (26,67,0,16,97) & (100,18,0,19,91) & (71,12,0,75,127) & (0,45,18,82,110) & mod 129 & 3354 & $\frac{13\times 129}{3354}$\\
  & (101,36,0,95,79) & (0,82,11,21,34) & (0,12,49,110,63) & (57,122,79,0,120) & & &\\
  & & & & & & &\\
  & (0,54,57,59,108) & (0,83,114,11,123) & (11,0,105,111,71) & (8,26,$\infty$,51,0) & & &\\
  & (3,0,112,70,69) & (0,124,41,101,74) & (0,24,32,46,122) & (40,52,36,0,56) & & &\\
  & & & & & & & \\
  & (1,32,33,0,116) & (0,26,50,94,53) & (0,39,68,77,35) & (0,69,39,75,92) & & & \\
  & (10,0,101,7,20) & (23,11,0,85,127) & (0,15,105,42,122) & (0,30,5,38,114) & & & \\
  & & & & & & & \\
  & (95,0,21,40,120) & & & & & & \\
  & (8,93,0,108,52) & & & & & & \\
  \hline
  131 & (0,70,1,17,59) & (0,118,2,34,9) & (68,0,4,18,105) & (8,36,5,0,79) & mod 131 & 3406 & $\frac{13\times 131}{3406}$\\
& (70,53,11,69,0) & (32,20,0,13,54) & (26,49,0,117,93)  &  (55,0,52,98,103) & & &\\
& & & & & & &\\
& (16,0,72,10,27) & (0,40,108,26,64) & (104,125,39,0,29) & (78,119,58,77,0) & & &\\
& (72,62,45,56,0) & (18,81,0,14,44) & (100,71,92,0,35) & (77,130,89,0,19) & & &\\
& & & & & & &\\
& (0,25,107,23,116) & (101,50,0,46,83) & (22,9,0,106,7) & (0,28,36,88,31) &  & &\\
& (129,23,47,0,38) & (127,76,46,0,94) & (0,124,13,112,90) & (128,80,52,0,85) & & &\\
& & & & & & &\\
& (123,21,92,0,57) & & &  & & &\\
&  (0,104,65,110,75) & & & &  & &\\
  \hline
\end{tabular}
\end{center}

\end{document}